\documentclass[11pt]{amsart}
\usepackage{amssymb}
\usepackage{amsfonts} 
\usepackage{amsmath}
\usepackage{mathrsfs}
\usepackage[margin=1.2in, centering]{geometry}
\numberwithin{equation}{section}

\newtheorem{theorem}{Theorem}[section]
\newtheorem{lemma}[theorem]{Lemma}

\def \mco {{\mathcal O}}
\def \mcp {{\mathcal P}}

\def \mcs {{\mathcal S}}

 \def \msf {{\mathscr F}}
 \def \mss {{\mathscr S}}
 \def \msl {{\mathscr L}}
 \def \msb {{\mathscr B}}

\def \mbc {{\mathbb C}}

\def \mbr {{\mathbb R}}
\def \mbs {{\mathbb S}}

\def \beqq {\begin{equation}}
\def \eeqq {\end{equation}}
\def \bpf {\begin{proof}}
\def \epf {\end{proof}}
\def \beq {\begin{equation*}}
\def \eeq {\end{equation*}}

\def \comp {\operatorname{comp}}
 
\def \im {\operatorname{Im}}

\def \phg {\operatorname{phg}}

\def \eps {\epsilon}   
\def \la {\lambda}   
       
\def \lap {\Delta}
\def \p {\partial}
\def \ha {\frac{1}{2}}

\begin{document} 
\title[]{Recovering asymptotics of potentials from the scattering of relativistic Schr\"odinger operators}
\author{Gunther Uhlmann}
\address{Gunther Uhlmann
\newline
\indent Department of Mathematics, University of Washington}
\email{gunther@math.washington.edu}
\author{Yiran Wang}
\address{Yiran Wang
\newline
\indent Department of Mathematics, Emory University}
\email{yiran.wang@emory.edu}
\begin{abstract}
We study the stationary scattering for $(-\lap)^\ha + V(x)$ on $\mbr^3$. For poly-homogeneous potentials decaying at infinity, we prove that the asymptotics of the potential can be recovered from the scattering matrix at a fixed energy. 
\end{abstract}
\date{\today}
\maketitle

\section{Introduction}
Let $(-\lap)^\ha$ be the fractional Laplacian on $\mbr^3$ defined by $(-\lap)^\ha f(x) = (|\xi| \widehat f(\xi))^\vee,$ 
where $\wedge, \vee$ denote the Fourier and inverse Fourier transform, respectively.  We consider the Hamiltonian 
\beqq\label{eq-ham}
H_V = (-\lap)^\ha + V(x), 
\eeqq 
where $V(x)$ is a  potential function. In the literature, \eqref{eq-ham} is called the relativistic Schr\"odinger operator or fractional Schr\"odinger operator. We refer to the introduction of  \cite{CMS} for a discussion of the  terminology and physical background. Recently, this operator was introduced in quantum optics as a model for the interaction between a quantized field and a collection of two-level atoms, see \cite{KrSc, LNOS}. 

For the time dependent problem, both the scattering and inverse scattering problems have been studied in the literature, see for example \cite{Ish, Jun}. Here, we are interested in the stationary scattering at a fixed energy. In this work, we assume that  $V$ is real-valued and satisfies 
\beqq\label{eq-pot}
|V(x)|\leq C \langle x \rangle^{-\sigma}, \quad \sigma >3
\eeqq
for $x\in \mbr^3$. The decay condition guarantees a well-defined scattering matrix for \eqref{eq-ham}.  It is known (see e.g.\ Theorem 4A of \cite{BeNe}) that $H_V$ is essentially self-adjoint  on $L^2(\mbr^3)$ with essential spectrum $\sigma_{ess}(H_V) = [0, \infty)$. The point spectrum $\sigma_{p}(H_V)$ on $[0, \infty)$ is discrete and can only accumulate at $0, \infty$. At present, the absence of embedded eigenvalues is only known under certain conditions on the potential, see \cite{ILS}. For $\la \in (0, \infty)\backslash \sigma_{p}(H_V)$, we consider solutions of $(H_V - \la) u = 0$ on  $\mbr^3$. For $u$ in certain weighted $L^2$ spaces, the equation holds in the sense of distributions. In Section \ref{sec-sca}, we show that $u$ has the following asymptotics as $|x|\rightarrow \infty$
\beqq\label{eq-uexpan}
u(|x|\theta)  =  e^{-i\la|x|} |x|^{-1} g(\theta) +  e^{i\la|x|} |x|^{-1} g'(\theta) + o_{L^\infty}(|x|^{-1}), 
\eeqq
where $\theta = x/|x|, x\in \mbr^3$. In fact, $g'$ is uniquely determined by $g$. Similarly to the potential scattering for the Laplacian, we define the scattering matrix $S_V(\la)$ to be   
\beqq\label{eq-scat}
S_V(\la)g = g'. 
\eeqq

We study the inverse problem of determining $V$ from $S_V(\la)$. For $m \in \mbr$, we say that $V \in S^m_{\phg}(\mbr^3)$ is polyhomogeneous of order $m$ if $V\in C^\infty(\mbr^3)$ and 
\beq
V(x) \sim \sum_{j = 0}^\infty V_j(x), 
\eeq
where $V_j$ are homogeneous of degree $-(m + j)$ for $|x|>1$, see Definition 18.1.5 of \cite{Ho3}. We set $S^{3+}_{\phg}(\mbr^3) = \bigcup_{m >3} S^m_{\phg}(\mbr^3)$ and $S^{-\infty}_{\phg}(\mbr^3) = \bigcap_{m\in \mbr}S^m_{\phg}(\mbr^3)$.  Our main result is  
\begin{theorem}\label{thm-main}
Let  $V_1, V_2 \in S^{3+}_{\phg}(\mbr^3)$ and $\la\in (0, \infty)\backslash (\sigma_p(H_{V_1})\cup \sigma_p(H_{V_2}))$. If $S_{V_1}(\la) = S_{V_2}(\la)$, then $V_1 - V_2\in S^{-\infty}_{\phg}(\mbr^3)$. 
\end{theorem} 
 
From the proof, we can see that the conclusion holds if the difference of the corresponding scattering amplitudes (i.e.\ Schwartz kernels of the scattering matrices) is smooth.  Also, our method works for potentials with more general asymptotic behaviors  as those in \cite{WeYa}.  Theorem \ref{thm-main} provides unique determination results for certain non-compactly supported potentials. For example, potentials that are rational functions with suitable decay at infinity can be uniquely determined from the scattering matrix of \eqref{eq-ham} at a fixed energy. 

We review some previous results and explain how we prove Theorem \ref{thm-main}. For the potential scattering for the Laplacian $-\lap$, recovering asymptotics of potentials has been studied extensively. We refer to the introduction of \cite{WeYa} for a summary. The general philosophy is that the asymptotics of potentials are encoded in the singularities of the scattering amplitude. For potentials with suitable behaviors at infinity, it can be shown that the scattering matrix is a Fourier Integral Operator (FIO) associated with the geodesic flow at infinity, see Theorem 1.1 of \cite{JoSa} and Theorem 3.2 of \cite{WeYa}. One can then recover the asymptotics of potentials from the principal symbol.  Typically, the proof relies on some construction of the Poisson operator as an oscillatory integral. In \cite{JoSa} and several subsequent works, this is achieved by using the calculus of Legendrian distributions in \cite{MeZw}. For the Euclidean scattering, a more direct construction was shown in \cite{WeYa}. For the scattering of \eqref{eq-ham}, we follow the same philosophy but difficulty arises from the fact that $(-\lap)^\ha$ is a non-local operator. Although $(-\lap)^\ha$ can be regarded as a pseudo-differential operator, the previously mentioned constructions of the Poisson operator seem difficult to carry out. Instead, we will compute the scattering amplitude directly from the generalized eigenfunction using the expansion of the resolvent, see \eqref{eq-expan}. In particular, for poly-homogeneous potentials, we show that  the leading order singularities of the scattering amplitude can be identified from the first term of expansion, which is sufficient for proving Theorem \ref{thm-main}. 

The paper is organized as follows. In Section \ref{sec-sca}, we discuss the definition of the scattering matrix for \eqref{eq-ham}. In Section \ref{sec-sing}, we analyze the leading order singularities in the scattering amplitude. Then we prove Theorem \ref{thm-main} in Section \ref{sec-pf}. Finally, in Section \ref{sec-bdy}, we prove a boundary pairing lemma that is used in the analysis.

\section{The scattering matrix}\label{sec-sca}  
\subsection{Free scattering}\label{subsec-free}
Consider $H_0 = (-\lap)^\ha$. Unlike the local operator $-\lap$, the symbol $|\xi|$ of $H_0$  is singular at $0$. As a result, even for $\phi \in \mss(\mbr^3)$,  $H_0 \phi$ may not belong to $\mss(\mbr^3)$.  The mapping properties of $H_0$ on weighted Sobolev spaces was studied in \cite{Ume1}. For $x\in \mbr^3$, 
$|x|$ denotes the Euclidean norm of $x$, and $\langle x \rangle = (1 + |x|^2)^\ha$. 
$D_j$ stands for $-i \p/\p x_j, j = 1, 2, 3$ and $\langle D\rangle = (1+ |D|^2)^\ha = (1 - \lap)^\ha$.  
For $k, s\in \mbr$, the weighted Sobolev space $H^{k, s}(\mbr^3)$ is defined by 
\beq
H^{k, s}(\mbr^3)  = \{f\in \mss'(\mbr^3)| \langle x\rangle^s\langle D\rangle^k f\in L^2(\mbr^3)\} 
\eeq
with the norm $\|f\|_{k, s} = \|\langle x\rangle^s \langle D\rangle^k f\|_{L^2}$. When $k = 0$, $H^{0, s}(\mbr^3)$ is the weighted $L^2$ space 
\beq
L^{2, s}(\mbr^3)  = \{f\in \mss'(\mbr^3)| \langle x\rangle^s  f\in L^2(\mbr^3)\} 
\eeq
with the norm $\|f\|_{s} = \|\langle x\rangle^s  f\|_{L^2}$. In Theorem 4.4 and 4.6 of \cite{Ume1}, it is proved that for $s<5/2$, $H_0$ maps $\mss(\mbr^3)$ continuously to $L^{2, s}(\mbr^3)$, and for $f\in L^{2, -s}(\mbr^3)$, $H_0 f\in \mss'(\mbr^3)$. The mapping properties on weighted Sobolev spaces are also obtained in \cite{Ume1}. 

Now for $w\in \mbs^2, \la>0$, let $\Phi_0(x; w, \la) = e^{i \la x\cdot w}$ be a plane wave. We know $\Phi_0\in L^{2, s}(\mbr^3)$ for any $s< -1$. Thus $H_0\Phi_0$ makes sense in  $\mss'(\mbr^3)$ and we have 
\beq
(H_0 - \la) \Phi_0 = 0 \text{ in } \mss'(\mbr^3),
\eeq
see Lemma 8.1 of \cite{Ume}. 
For $h\in C^\infty(\mbs^2)$, we define 
\beq
\Phi_0(\la)h(x) = \int_{\mbs^2} e^{i \la x\cdot w} h(w)dw.
\eeq
We can apply the stationary phase lemma to get a complete asymptotic expansion for $|x|\rightarrow \infty$
\beqq\label{eq-phiasymp}
\begin{gathered}
\Phi_0(\la)h(|x|\theta) \sim e^{-i\la |x|} (\la |x|)^{-1} (2\pi i) \sum_{j =0}^\infty |x|^{-j} h_j^-(\theta)
+e^{i\la |x|} (\la|x|)^{-1} (-2\pi i) \sum_{j = 0}^\infty |x|^{-j} h_j^+(\theta). 
\end{gathered}
\eeqq
We refer to Section 1.3 of \cite{Mel} for the details.  In particular, $h_j^\pm(\theta), j\geq 1$ are all given by polynomials in the Laplacian on $\mbs^2$ applied to $h(\pm\theta)$.  Here, we recall that $\sim$ means that for any $N>0$ integer and $|x|>1$, 
\beq
\begin{gathered}
| \Phi_0(\la)h(|x|\theta)-  e^{-i\la |x|} (\la |x|)^{-1} (2\pi i) \sum_{j =0}^N |x|^{-j} h_j^-(\theta) 
- e^{i\la |x|} (\la|x|)^{-1} (-2\pi i) \sum_{j = 0}^N|x|^{-j} h_j^+(\theta)| \leq C_N |x|^{-N-2}
\end{gathered}
\eeq
for some $C_N>0$. To show that the expansion \eqref{eq-phiasymp} is unique in a proper sense, we need the following boundary pairing lemma. 
\begin{lemma}\label{lm-bp2}
Suppose  $\la>0$, $((-\lap)^\ha - \la) u_\pm \in L^{2, s}(\mbr^3), s\geq 1$ and
\beqq\label{eq-upm}
\begin{gathered}
u_+ = e^{-i\la |x|} |x|^{-1}g_{+-}(\theta) + e^{i \la |x|}|x|^{-1}g_{++}(\theta) + L^2(\mbr^3),\\
u_- = e^{-i\la |x|} |x|^{-1}g_{--}(\theta) + e^{i \la |x|}|x|^{-1}g_{-+}(\theta) + L^2(\mbr^3), 
\end{gathered}
\eeqq
where $g_{\pm\pm}\in C^\infty(\mbs^2)$. Then
\beq
\begin{gathered}
\langle u_+, ((-\lap)^\ha  - \la)u_-\rangle -  \langle ((-\lap)^\ha  - \la)u_+, u_-\rangle 
= i(\langle g_{++}, g_{-+}  \rangle - \langle g_{+-}, g_{--} \rangle).
\end{gathered}
\eeq    
\end{lemma} 

We postpone the proof to Section \ref{sec-bdy}. Now we have 
\begin{lemma}\label{lm-gen0}
For $\la \in(0, \infty)$ and $h\in C^\infty(\mbs^2)$, there is a unique solution $u\in \msl \doteq \bigcup_{s<5/2}L^{2, -s}(\mbr^3)$ to $(H_0 - \la)u = 0$ in $\mss'(\mbr^3)$   such that as $|x|\rightarrow \infty$, 
\beqq\label{eq-uasy} 
u(|x|\theta) = e^{-i\la |x|} |x|^{-1}h(\theta) + e^{i \la |x|}|x|^{-1} h'(\theta) + O_{L^\infty}(|x|^{-2}), 
\eeqq
where $h' \in C^\infty(\mbs^{2})$ and $h'(\theta) = -h(-\theta).$ We define the scattering matrix $S_0(\la): C^\infty(\mbs^2)\rightarrow C^\infty(\mbs^2)$ by $S_0(\la)h = h'$. 
\end{lemma} 
\bpf
For $h\in C^\infty(\mbs^2)$, we let $h^*(\theta) = h(-\theta)$ for $\theta\in \mbs^2.$ We set $u = \Phi_0(\la)h^*$ and observe that $u\in L^{2, -3/2}(\mbr^3)$ so $u\in \msl,$ and $(H_0 - \la)u = 0$ in $\mss'(\mbr^3)$. The asymptotic expansion is given by \eqref{eq-phiasymp}. So the existence part is done. For the uniqueness, let $\tilde u$ be another solution of the form \eqref{eq-uasy}, that is 
\beq
\tilde u(|x|\theta) = e^{-i\la |x|} |x|^{-1}h(\theta) + e^{i \la |x|}|x|^{-1} \tilde h'(\theta) + O_{L^\infty}(|x|^{-2}), 
\eeq
where $\tilde h'\in C^\infty(\mbs^2)$. Then the difference 
\beq
w(|x|\theta)\doteq u(|x|\theta) - \tilde u(|x|\theta) = e^{i \la |x|}|x|^{-1}h''(\theta) + O_{L^\infty}(|x|^{-2}), 
\eeq
where $h'' = h' - \tilde h'\in C^\infty(\mbs^2)$. Note that $((-\lap)^\ha -\la)w = 0$. We can apply Lemma \ref{lm-bp2} to get 
\beq
\begin{gathered}
 0   = \langle w, ((-\lap)^\ha  - \la)w\rangle -  \langle ((-\lap)^\ha  - \la)w, w\rangle = i \langle h'', h''\rangle. 
\end{gathered}
\eeq  
We conclude that $h''=0$.  Then $u - \tilde u\in L^2(\mbr^3)$.  Because $\la$ is not an eigenvalue of $H_0$, $u = \tilde u.$
\epf

\subsection{Potential scattering}  \label{subsec-pot}
It is known (see e.g.\ \cite{BeNe}) that $H_0$ is a self-adjoint operator on $L^2(\mbr^3)$ with domain $H^1(\mbr^3)$. The spectrum of $H_0$ is $\sigma_{ess}(H_0) = [0, \infty)$. 
For $V$ satisfying \eqref{eq-pot}, $V$ is relatively compact with respect to $H_0$ (in fact, $\sigma>1$ is sufficient, see \cite{Ume}). It follows that $H_V$, as a self-adjoint operator on $L^2(\mbr^3)$ with domain $H^1(\mbr^3)$, has essential spectrum $\sigma_{ess}(H_V) = [0, \infty)$. Also, the point spectrum on $(0, \infty)$ is a discrete set with possible accumulation points at $0, \infty$, see \cite{BeNe}.  For $z \in \mbc\backslash(\sigma_{ess}(H_V)\cup \sigma_{p}(H_V))$, we let $R_V(z) = (H_V - z)^{-1}$ be the resolvent.  The limiting absorption principle was proved in \cite{BeNe}. For $\la \in (0, \infty)\backslash \sigma_{p}(H_V)$ and $s>\ha$, there exists the limit
\beqq\label{eq-lim0}
R_V^\pm(\la) = \lim_{\mu\rightarrow 0+} R_V(\la \pm i\mu) \text{ in } \msb(L^{2, s}, H^{1, -s}), 
\eeqq  
where $\msb(A, B)$ denotes the set of bounded operators from $A$ to $B$.  In particular, we can define the operator-valued functions $R_V^\pm(z)$ by 
\beqq\label{eq-lim}
R_V^\pm(z) = \begin{cases}
&R_V(z), \quad z\in \mbc^\pm, \\
& R_V^\pm(z), \quad z = \la >0. 
\end{cases}
\eeqq
Then $R_V^\pm(z)$ are $\msb(L^{2, s}, H^{1, -s})$ valued continuous functions.

We consider the construction of generalized eigenfunctions for $H_V$. For $w\in \mbs^2, \la>0$, let $\Phi_0(x; w, \la) = e^{i \la x\cdot w}$ be a plane wave. 
In \cite[Section 8]{Ume}, generalized eigenfunctions for $H_V$ are defined by
\beqq\label{eq-resV}
\Phi^\pm_V = \Phi_0 - R_V^\pm(\la) (V\Phi_0)
\eeqq
which belongs to $L^{2, -s}(\mbr^3), s>3/2$. Here, we changed the sign convention in \cite{Ume}. For $\la \in (0, \infty)\backslash \sigma_{p}(H_V)$, both $\Phi^\pm_V$ satisfy 
\beq
(H_V - \la) \Phi^\pm_V = 0 \text{ in } \mss'(\mbr^3). 
\eeq 
Furthermore, in Theorem 10.2 of \cite{Ume}, it is proved that 
\beqq\label{eq-phiasymp1}
\Phi^\pm_V(|x|\theta; w, \la) = \Phi_0(|x|\theta; w, \la) + e^{\pm i \la |x|}{|x|^{-1}} f_V^{\pm}(\theta, w, \la) + \gamma^\pm_V(x; w, \la), 
\eeqq
where 
\beqq\label{eq-f}
f_V^\pm(\theta, w, \la) = -\frac{\la}{2\pi} \int_{\mbr^3} e^{\pm i\la \theta \cdot y} V(y) \Phi_V^\pm(y;  w, \la) dy
\eeqq
and 
\beq
|\gamma_V^\pm(x; w, \la)|\leq C \begin{cases} 
&|x|^{-(\sigma-1)/2}, \quad 3< \sigma < 5\\
&|x|^{-2} \log(1+ |x|), \quad \sigma = 5\\
&|x|^{-2}, \quad \sigma >5
\end{cases}
\eeq 
Note that $\gamma_V^\pm \in o_{L^\infty}(|x|^{-1})$ as $|x|\rightarrow \infty.$ For each $w\in \mbs^2,$ $\Phi^\pm_V(x; w, \la)$ are continuous functions in $x\in \mbr^3$ and $\la\in (0, \infty)\backslash \sigma_p(H_V)$, see Theorem 9.1 of \cite{Ume}. We see from \eqref{eq-f} that $f_V^\pm$ are continuous in $\theta, w.$  

We are ready to define the scattering matrix for $H_V.$
\begin{lemma}\label{lm-gen}
Suppose $V$ satisfies \eqref{eq-pot}. For each $\la \in (0, \infty)\backslash \sigma_{p}(H_V)$ and $h\in C^\infty(\mbs^2)$, there is a unique solution in $\msl$ to $(H_V - \la)u = 0$ on $\mss'(\mbr^3)$ with the following asymptotics as $|x|\rightarrow \infty$ 
\beqq\label{eq-ufull}  
u(|x|\theta) =  e^{-i\la |x|} |x|^{-1} h(\theta) + e^{i \la |x|}|x|^{-1}h'(\theta) +  o_{L^\infty}(|x|^{-1}), 
\eeqq
where  $h' \in L^\infty(\mbs^{2})$ is determined by $h$ via  
\beq
h'(\theta) = -h(-\theta) + \int_{\mbs^2} f_V^{+}(\theta, w, \la) h(-w)dw. 
\eeq
The scattering matrix is defined by $S_V(\la): C^\infty(\mbs^2)\rightarrow L^\infty(\mbs^2)$, $S_V(\la) h = h'.$ We call $f_V^+$ the scattering amplitude. 
\end{lemma}

\bpf 
For $h\in C^\infty(\mbs^2)$, we let $h^*(\theta) = h(-\theta),$ for $\theta\in \mbs^2.$ We set $u = \Phi_V^+(\la)h^*$. Using \eqref{eq-phiasymp1}, we know that $u\in L^{2, -3/2}(\mbr^3)$ so $u\in \msl$. Also, $(H_V - \la)u = 0$ in $\mss'(\mbr^3)$. Using \eqref{eq-phiasymp}, we get that 
\beq
u(x) = \int_{\mbs^2}\Phi_0(|x|\theta; w, \la)h(-w)dw + \int_{\mbs^2} e^{i \la |x|}{|x|^{-1}} f_V^{+}(\theta, w, \la) h(-w)dw + \int_{\mbs^2} \gamma_V^+(x; w, \la) h(-w)dw. 
\eeq
We can apply the stationary phase lemma for the first integral to get the asymptotic expansion. The existence part is done. 

For the uniqueness part, let $\tilde u$ be another solution of the form \eqref{eq-ufull}, that is 
\beq
\tilde u(|x|\theta) = e^{-i\la |x|} |x|^{-1}h(\theta) + e^{i \la |x|}|x|^{-1} \tilde h'(\theta) + o_{L^\infty}(|x|^{-1}), 
\eeq
where $\tilde h'\in L^2(\mbs^2)$. Then the difference 
\beq
w(|x|\theta) \doteq u(|x|\theta) - \tilde u(|x|\theta) = e^{i \la |x|}|x|^{-1}h''(\theta) + o_{L^\infty}(|x|^{-1}), 
\eeq
where $h'' = h' - \tilde h'\in L^2(\mbs^2)$. Note that $((-\lap)^\ha + V-\la)w = 0$ so $((-\lap)^\ha - \la)w = -Vw \in L^{2, 1}(\mbr^3)$. We can apply Lemma \ref{lm-bp2} to get 
\beq
\begin{gathered}
 0 =  \langle w, ((-\lap)^\ha + V-\la) w \rangle -  \langle ((-\lap)^\ha  + V - \la)w, w\rangle \\
 = \langle w, ((-\lap)^\ha  - \la)w\rangle -  \langle ((-\lap)^\ha  - \la)w, w\rangle = i \langle h'', h''\rangle. 
\end{gathered}
\eeq  
We conclude that $h''=0$ so $u - \tilde u\in L^2(\mbr^3)$.  Because $\la$ is not an eigenvalue of $H_V$, $u = \tilde u.$
\epf

\section{Singularities of the scattering amplitude}\label{sec-sing}
We analyze the singularities in $f_V^+.$ In view of \eqref{eq-phiasymp1} and \eqref{eq-resV}, it suffices to examine the asymptotics of $w = R_V^+(\la)(V\Phi_0)$. For $\im \la >0$, we have the resolvent identity 
\beq
R_V(\la) = R_0(\la) - R_0(\la) VR_V(\la). 
\eeq
Using the limiting absorption principle \eqref{eq-lim0} and \eqref{eq-lim}, we get 
\beqq\label{eq-born0}
w = R_0^+(\la)(V\Phi_0) - R_0^+(\la)(Vw)
\eeqq
for $\la\in (0, \infty)\backslash \sigma_p(H_V)$. Here, we note that  $V\Phi_0\in L^{2, 1}(\mbr^3)$ hence we know a priori $w\in H^{1, -1}(\mbr^3)$ by \eqref{eq-lim0}. Thus $Vw\in L^{2, 1}(\mbr^3)$ by the decay property of $V$ so \eqref{eq-born0} makes sense. The formula \eqref{eq-born0} is sufficient for our analysis. But we remark that we can use the formula iteratively to formally obtain the expansion 
\beqq\label{eq-expan}
w = R_0^+(\la)(V\Phi_0) + \sum_{j = 1}^\infty (-R_0^+(\la) V)^j \Phi_0. 
\eeqq
In the following, we assume that $V$ is of the form  
\beqq\label{eq-pot1}
V(x) = V_0(x) + V_1(x), 
\eeqq
where for some $m>3$ and $\beta>0$ small, we have 
\begin{enumerate}
\item $V_0\in C^\infty(\mbr^3)$ and $V_0(x) = |x|^{-m} v_0(x/|x|)$ for $|x|>1$. 
\item $V_1\in L^\infty(\mbr^3)$ and $V_1 \in O_{L^\infty}(|x|^{-m-\beta})$ for $|x|\rightarrow \infty$. 
\end{enumerate}
Note that $V$ is not necessarily poly-homogeneous.  We decompose \eqref{eq-born0} as 
\beqq\label{eq-born1}
\begin{split}
w &= R_0^+(\la)(V_0 \Phi_0) + R_0^+(\la)(V_1 \Phi_0) - R_0^+(\la)(Vw) \\
 & \doteq w_0 + w_1 + w_2.
 \end{split}
\eeqq
We will show in Lemma \ref{lm-uasymp} and \ref{lm-uasymp1} that $w_0$ contains the leading order singularity of the scattering amplitude, and $w_1, w_2$ are more regular. We remark that for the purpose of solving the inverse problem,  it is sufficient to find the characterization of singularities in the Sobolev scale. It is not necessary to show the fine microlocal structure of the singularities, although this can be achieved at least for the singularities in $w_0$ by our approach. 

Before getting to  the lemmas, we discuss the X-ray transform of $V_0$. For $\theta\in \mbs^2$, we let $\theta^\perp = \{x\in \mbr^3: x\cdot \theta = 0\}$. Then all straight lines in $\mbr^3$ can be parametrized as $x + t\theta, t\in \mbr$, where $\theta\in \mbs^2, x\in \theta^\perp$. We define the X-ray transform of $V_0$ as 
\beqq\label{eq-xray}
X V_0(x, \theta) = \int_\mbr V_0(x + t\theta)dt, \quad \theta \in \mbs^2, x\in \theta^\perp. 
\eeqq
For fixed $\theta$, we consider the asymptotic behavior of $XV_0(x, \theta)$ as $|x|\rightarrow \infty.$ Without loss of generality, we take $\theta = (0, 0, 1)$ so $\theta^\perp = \{(x', 0): x' =(x_1, x_2) \in \mbr^2\}$. Let $w \in \mbs^2$ and $w\cdot \theta = 0$. We compute the asymptotics of $XV_0(rw, \theta)$ as $r\rightarrow \infty$.  We can choose the coordinates so that $w  = (1, 0, 0)$. For $r>1$, we have 
\beq
X V_0(rw, \theta) = \int_\mbr V_0(r w + t\theta)dt = \int_\mbr \frac{1}{(r^2 + t^2)^{m/2}}v_0(\frac{r}{(r^2 + t^2)^\ha}, 0, \frac{t}{(r^2 + t^2)^\ha}) dt. 
\eeq
Making a change of variable $t = rs$,  we get 
\beq
X V_0(rw, \theta)  =  r^{-m+1}\int_\mbr \frac{1}{(1 + s^2)^{m/2}}v_0(\frac{1}{(1 + s^2)^\ha}, 0, \frac{s}{(1 + s^2)^\ha})ds.
\eeq
Then we make another change of variable $\sin \alpha = 1/(1+s^2)^\ha, \cos \alpha = s/(1+s^2)^\ha, \alpha \in (0, \pi)$ to get 
\beqq\label{eq-rays2}
X V_0(rw, \theta)  =  r^{-m+1}\int_{0}^{\pi} (\sin \alpha)^{m-2}v_0(\sin \alpha, 0, \cos \alpha) d\alpha. 
\eeqq
We denote the above integral by $I(\theta, w)$ and write 
\beqq\label{eq-I}
XV_0(rw, \theta) = r^{-m+1}I(\theta, w), 
\eeqq
so $XV_0(x, \theta)$ is homogenous of degree $-m+1$ in $x$ for $x\in \theta^\perp$ and $|x|>1$. Furthermore, the integral $I$ has the following interpretation. We identify $\theta^\perp$ with $T_{-\theta}\mbs^2$ and  $w$ as a tangent vector in $T_{-\theta}\mbs^2$. Then  
\beq
\gamma_{\theta, w}(\alpha) = (\sin \alpha, 0, \cos \alpha), \quad \alpha\in (\pi, 0)
\eeq
is the unit speed geodesic on $\mbs^2$ from $-\theta$ in direction $w$ to $\theta$. Thus $I(\theta, w)$ is really a weighted geodesic ray transform of $v_0$ on $\mbs^2$. In \cite{JoSa}, such integrals appeared in the principle symbol of the scattering matrix, and the authors recovered $v_0$ by considering the geodesic ray transform on $\mbs^2$. We believe that the same can be done here. However, because of \eqref{eq-I}, knowing $I(\theta, w)$ is equivalent to knowing $X V_0$, and the latter is easier to invert. We remark that this consideration was used in \cite{WeYa} for the potential scattering of $-\lap$. 

\begin{lemma}\label{lm-uasymp}
Let $\la \in (0, \infty)\backslash \sigma_{p}(H_V)$. Then as $|x| \rightarrow \infty,$ we have
\beqq\label{eq-RV}
w_0(|x|\theta) =   |x|^{-1} e^{i \la |x|} f_0(\la, \theta, w)  + L^2(\mbr^3), \quad \theta, w\in \mbs^2, 
\eeqq 
where $f_0(\la, \theta, w) = (-i) (2\pi)^{-2} \la  \hat V_0(\la(\theta -  w))$ is smooth away from $\theta = w$. Moreover, for  fixed  $\theta\in \mbs^2,$ $f_0$ regarded as a function of $w\in \mbs^2$ belongs to $H^{\kappa}(\mbs^2)$ for any $\kappa < m-2$. If $XV_0(x, \theta)$ is not identically zero for $x\in \theta^\perp, |x|>1$, then $f_0$ is not in $H^{m -2}(\mbs^2)$. 
\end{lemma}

\bpf
We first compute the asymptotics of $w_0(x)$ as $|x| \rightarrow \infty.$ Because $V_0 \Phi_0\in L^{2, s}(\mbr^3)$ for $s< 2$, we know from \eqref{eq-lim0} that $w_0\in H^{1, -s}(\mbr^3)$. For $\eps>0$, we let $\phi_\eps(x) = e^{-\eps|x|}, x\in \mbr^3$. Then we set $V_\eps = V_0 \phi_\eps$. We will compute the asymptotics of $w_{0, \eps}(x) = R_0^+(\la)(V_\eps \Phi_0)(x)$ as $|x|\rightarrow\infty$ and in the end, we will show that $w_{0, \eps} \rightarrow w_0$ in the sense of distributions as $\eps\rightarrow 0+.$
 
By using the Fourier transform, we find that 
 \beqq\label{eq-u1}
\begin{gathered}
w_{0, \eps}(x) = (2\pi)^{-3} \int_{\mbr^3} e^{i x \cdot \xi} \frac{1}{|\xi| - \la}   \hat V_\eps(\xi - \la w)  d\xi. 
\end{gathered}
\eeqq 
We compute the integration in polar coordinates. Let $x= r\theta, r>0, \theta\in \mbs^2$ and $\xi = \rho \eta, \rho>0, \eta\in \mbs^2$. We have 
\beq
\begin{gathered}
w_{0, \eps}(r\theta) = (2\pi)^{-3} \int_{\mbs^2}\int_0^\infty e^{i r\theta \cdot \rho \eta} \frac{1}{\rho - \la} \hat V_\eps(\rho\eta - \la w) \rho^2  d\rho d\eta. 
\end{gathered}
\eeq
Let $\phi$ be a smooth cut-off function on $(0, \infty)$ so that $\phi(\rho) = 1$ for $\rho\in (0, \la/2)$. We see that 
\beq
\begin{gathered}
w_{0, \eps}(r\theta) = (2\pi)^{-3} \int_{\mbs^2}\int_0^\infty e^{i r\theta \cdot \rho \eta} \frac{1}{\rho - \la} \hat V_\eps(\rho\eta - \la w) \rho^2 (1 - \phi(\rho))  d\rho d\eta + L^2(\mbr^3). 
\end{gathered}
\eeq
Note that because $V_\eps$ decays exponentially fast at infinity, $\hat V_\eps(\zeta)$ is analytic for $|\im\zeta| < \eps$, see Theorem IX.3 of \cite{ReSi2}.  The integral in $\rho$ can be computed by a contour deformation. We consider a cut-off function $\psi_0$ on $\mbs^2$ such that $\psi_0(w)$ is identically equal to $1$ in a neighborhood of $w\cdot \theta = 0$, and $\psi_0(w)$ vanishing identically in neighborhoods of both $w = \pm \theta$. Let $\psi_\pm$ be supported in $\pm \theta \cdot w \geq 0$ such that $\psi_0, \psi_\pm$ form a partition of unity with $\psi_0 + \psi_++\psi_- = 1$ on $\mbs^2$. Then we write 
\beq
\begin{gathered}
w_{0, \eps}(x) = v_{0, \eps}(x) + v_{+, \eps}(x) + v_{-, \eps}(x) + L^2(\mbr^3), \text{ where }\\
v_{t, \eps}(x) =  (2\pi)^{-3}  \int_{\mbs^2} \int_0^\infty  e^{i \rho r\theta \cdot \eta} \frac{1 - \psi(\rho \eta)}{\rho - \la} \psi_t(\eta)\hat V_\eps(\rho \eta - \la w) (1 - \phi(\rho))\rho^2 d\rho d\eta.   
\end{gathered}
\eeq
Since we are considering $R_0^+(\la)$, we analyze the above integral for $\im\la \geq0$. Of the three terms, $v_{0, \eps}, v_{-, \eps}\in L^2(\mbr^3)$ and this is uniform for $\im \la \geq 0$ and $\eps\geq 0$. For $v_{+, \eps}$, we can apply a contour deformation (see \cite[Proposition 1.1]{Mel} for the details) to get 
\beq
v_{+, \eps}(r\theta) = (2\pi)^{-3}  \int_{\mbs^2} e^{i \la r \theta \cdot \eta}  \psi_+(\eta) \hat V_\eps(\la \eta - \la w) \la^2  d\eta  + L^2(\mbr^3). 
\eeq
Now we can use the stationary phase lemma to compute the integration in $w$ to get for $r\rightarrow \infty$ that 
\beq
\begin{gathered}
v_{+, \eps}(r \theta)  = (-2\pi i) (2\pi)^{-3} \la r^{-1} e^{i \la r}  \hat V_\eps(\la\theta -  \la w) + L^2(\mbr^3), 
\end{gathered}
\eeq
in which the $L^2$ term is uniform for $\eps \geq 0.$ As $\eps\rightarrow 0+,$ we conclude that $v_{+, \eps}$ converges in the sense of distributions to 
\beq
\begin{gathered}
v_+(r \theta)  =  (-i)(2\pi)^{-2} \la r^{-1} e^{i \la r}  \hat V_0(\la\theta -  \la w) + L^2(\mbr^3).
\end{gathered}
\eeq

Next, for fixed $\theta\in \mbs^2$, we consider $U(w) \doteq \hat V_0(\la(\theta- w))$ as a function of $w\in \mbs^2$.  According to \eqref{eq-pot1}, we know that the Fourier transform $\hat V_0(\xi), \xi\in \mbr^3$ is in $H^{\kappa}(\mbr^3)$ for any $\kappa < m-3/2$.  In fact, because $V_0(x)$ is homogeneous of order $-m$, we know that $\hat V_0(\xi)$ is only singular at $\xi = 0$ (and the singularity is conormal). This implies that when restricted to $\mbs^2$, $U(w)$ is only singular at $w  = \theta.$ 

Finally, we consider the Sobolev regularity of $U(w)$ near $\theta$. Without loss of generality, we take $\theta = (0, 0, 1)\in \mbs^2$ and consider $w$ in a small neighborhood of $\theta$. For $\eps>0$, let 
\beqq\label{eq-o}
\mco = \{(w_1, w_2, (1 - w_1^2 - w_2^2)^\ha): |w_1|< \eps, |w_2|< \eps\}. 
\eeqq
For $w\in \mco$, we know that $\la(\theta - w)\in \mcs$ given by 
\beq 
 \mcs = \{-\la (w_1, w_2, (1- w_1^2 - w_2^2)^\ha - 1): |w_1|<\eps, |w_2|<\eps\}, 
\eeq
which is a smooth co-dimension one submanifold of $\mbr^3.$ Thus by the trace theorem,  $U|_\mco \in H^{\kappa }(\mco)$ for $\kappa < m-2$.

It remains to show that $U|_{\mco}$ is not in $H^{m-2}(\mco)$ if $XV_0$ is not identically zero. Using local coordinates in \eqref{eq-o}, we have  
\beq
\begin{gathered}
U(w_1, w_2)  = \int_{\mbr^3} e^{-i\la (w_1 x_1 + w_2 x_2 + x_3 ((1- w_1^2 - w_2^2)^\ha - 1)) } V_0(x_1, x_2, x_3)dx. 
 \end{gathered}
\eeq 
Write $w' = (w_1, w_2), x' = (x_1, x_2)$. We decompose $U = U_0 + U_1$ where 
\beq 
\begin{split}
U_0(w')  &= \int_{\mbr^2}\int_\mbr e^{-i\la x'\cdot w'} V_0(x', x_3)dx'dx_3  \\
&= \int_{\mbr^2} \int_{\mbr} e^{-i\la x'\cdot w'} (1 -\chi(x')) V_0(x', x_3)dx'dx_3 +  \int_{\mbr^2}\int_{\mbr} e^{-i\la x'\cdot w'} \chi(x') V_0(x', x_3)dx'dx_3 \\
&\doteq U_{0, 1}(w') + U_{0, 2}(w'). 
\end{split}
\eeq
Here, $\chi$ is a smooth cut-off function on $\mbr^2$ such that $\chi(x') = 1$ for $|x'|<1$ and $\chi(x') = 0$ for $|x'| > 2$. 
We consider $U_0, U_{0, 1}, U_{0, 2}$ as functions defined on $\mbr^2$. Note that $U_{0, 2}$ is Schwartz on $\mbr^2$. For $U_{0, 1},$ we compute  that 
\beqq\label{eq-U0}
\begin{split}
U_{0, 1}(w') &= \int_{\mbr^2} e^{-i\la w'\cdot x'} (1- \chi(x')) XV_0(x', \theta)dx'  \\
&=  \int_{\mbr^2} e^{-i\la w'\cdot x'} (1 - \chi(x')) \frac{1}{|x'|^{m-1}} I(\theta, (\frac{x'}{|x'|}, 0)) dx'
\end{split}
\eeqq
using \eqref{eq-I}. If $XV_0(x, \theta)$ is not identically zero, then according to \eqref{eq-I},  $I(\theta, (x'/|x'|, 0))$ is not identically zero.   The Fourier transform $\widehat{ U_{0, 1}}(x') = C (1 - \chi(x') ) I(\theta, (x'/|x'|, 0)) |x'|^{-m+1}$ where $C\neq 0$. We see that   $U_{0, 1}$ does not belong to $H^{m-2}(\mbr^2)$.  Note that $U_{0, 1}(w')$ for $|w'|>1$ is smooth and rapidly decaying at infinity. Thus, we conclude that $U_{0, 1}$ as well as $U_0$ do not belong to $H^{m-2}(\mco).$ 

For $U_1 = U - U_0$, we compute that 
\beq
\begin{split}
U_1(w')  &= \int_\mbr  \int_{\mbr^2} e^{-i\la w'\cdot x'}(e^{-i\la x_3 ((1- |w'|^2)^\ha - 1))} -1) V_0(x', x_3)dx'dx_3\\
&=  \int_\mbr  \int_{\mbr^2}  \frac{1}{-\la^2 |w'|^2} \lap_{x'} e^{-i\la w'\cdot x'}(e^{-i\la x_3 ((1- |w'|^2)^\ha - 1))} -1) V_0(x', x_3)dx'dx_3  \\
& = \int_\mbr   \int_{\mbr^2}    \frac{1}{\la^2 |w'|^2} e^{-i\la w'\cdot x'}(e^{-i\la x_3 ((1- |w'|^2)^\ha - 1))} -1) \lap_{x'} V_0(x', x_3)dx'dx_3. 
\end{split}
\eeq
Here, $\lap_{x'}$ denotes the Laplacian in the $x'$ variables. We estimate the integral in $x_3$
\beq
\begin{gathered}
| \frac{1}{\la^2 |w'|^2}  \int_{\mbr} (e^{-i\la x_3 ((1- |w'|^2)^\ha - 1))} -1) V_0(x', x_3) dx_3| \leq C   \int_{\mbr} |x_3 \lap_{x'} V_0(x', x_3)| dx_3 \\
\leq C  \int \frac{|x_3|}{(|x'|^2 + x_3^2)^{(m+2)/2}} dx_3\leq C \frac{1}{|x'|^{m}}. 
\end{gathered}
\eeq
Therefore, 
\beq
U_1(w') =  \int_{\mbr^2} e^{-i\la w'\cdot x'} \tilde V(x')dx', \text{ where } |\tilde V(x')|\leq C |x'|^{-m} \text{ for $|x'|$ large}.    
\eeq
This implies $U_1\in H^{m-2}(\mbr^2).$ So we proved that if $XV_0(x, \theta)$ is not identically zero, then $U$ does not belong to $H^{m-2}(\mco).$ This completes the proof of the lemma. 
\epf

\begin{lemma}\label{lm-uasymp1}
Let $\la \in (0, \infty)\backslash \sigma_{p}(H_V)$. Then as $|x| \rightarrow \infty,$ we have
\beqq\label{eq-RV1}
\begin{gathered}
w_j(|x|\theta)  =  |x|^{-1} e^{i \la |x|} f_j(\la, \theta, w)  + L^2(\mbr^3), \quad j = 1, 2,\quad \theta, w\in \mbs^2, 
\end{gathered}
\eeqq 
where  for fixed  $\theta\in \mbs^2,$ $f_j$ regarded as a function of $w\in \mbs^2$ belongs to $H^{m-2}(\mbs^2)$.  
\end{lemma}
\bpf
First, consider $w_1 =  R_0^+(\la)(V_1 \Phi_0)$ where $V_1\in O_{L^\infty}(|x|^{-m-\beta})$. Because $m>3$, we have $V_1\in L^2(\mbr^3)$. We can carry out the same calculation in Lemma \ref{lm-uasymp} to get 
\beq 
w_1(r\theta)  =  r^{-1} e^{i \la r} f_1(\la, \theta, w)  + L^2(\mbr^3), 
\eeq
where $f_1(\la, \theta, w) = (-i) (2\pi)^{-2} \la  \hat V_1(\la(\theta -  w))$. As in the proof of Lemma \ref{lm-uasymp}, we derive that $\hat V_1(\xi)\in H^{\kappa}(\mbr^3)$ for  $\kappa < m +\beta - 3/2$. Thus, for fixed $\theta$, $\hat V_1(\la(\theta-w))\in H^{\kappa}(\mbs^2), \kappa< m + \beta -2$. In particular, this implies that $f_1\in H^{m-2}(\mbs^2)$ because $\beta>0.$

Second, consider $w_2 = -R_0^+(\la)(V w)$. We note that $w(x) = O_{L^\infty}(|x|^{-1})$ so $Vw \in O_{L^\infty}(|x|^{-m-1})$. The same analysis for $w_1$ yields the conclusion of the lemma. This completes the proof. 
\epf

\section{Proof of Theorem \ref{thm-main}}\label{sec-pf}
Let  $V^{(j)} \in S^{3+}_{\phg}(\mbr^3), j = 1, 2$ and $\la\in (0, \infty)\backslash (\sigma_p(H_{V^{(1)}})\cup \sigma_p(H_{V^{(2)}}))$. Consider the corresponding generalized eigenfunctions $\Phi_{V^{(j)}}^+, j = 1, 2$ satisfying \eqref{eq-resV}. As in Section \ref{sec-sing}, we analyze the asymptotic of $w^{(j)} = R_{V^{(j)}}^+(\la)(V^{(j)}\Phi_0), j = 1, 2$. We have 
\beq 
((-\lap)^\ha - \la)w^{(j)} = V^{(j)}\Phi_0 - V^{(j)}w^{(j)} 
\eeq 
so that 
\beq
w^{(j)} = R_0^+(\la)(V^{(j)}\Phi_0) - R_0^+(\la)(V^{(j)}w^{(j)}). 
\eeq
Therefore, 
\beqq\label{eq-dif}
\begin{gathered}
w^{(1)} - w^{(2)} = R_0^+(\la)((V^{(1)} - V^{(2)})\Phi_0) - R_0^+(\la)((V^{(1)} - V^{(2)})w^{(1)}) \\
- R_0^+(\la)( V^{(2)}(w^{(1)} - w^{(2)})). 
\end{gathered}
\eeqq
For the left hand side of \eqref{eq-dif}, we have 
\beq
w^{(1)}(|x|\theta) - w^{(2)}(|x|\theta) 
=  o_{L^\infty}(|x|^{-1}) 
\eeq
using the assumption of Theorem \ref{thm-main} that $S_{V_1}(\la) = S_{V_2}(\la)$ and Lemma \ref{lm-gen}. Therefore, 
\beqq\label{eq-mainid}
\begin{gathered}
R_0^+(\la)((V^{(1)} - V^{(2)})\Phi_0) - R_0^+(\la)((V^{(1)} - V^{(2)})w^{(1)}) \\
- R_0^+(\la)( V^{(2)}(w^{(1)} - w^{(2)})) =   o_{L^\infty}(|x|^{-1}). 
\end{gathered}
\eeqq
We analyze the asymptotics of each terms. We write $V^{(j)}\in S^{3+}_{\phg}(\mbr^3), j = 1, 2$ in the form 
\beq
V^{(j)}(x) = V_0^{(j)}(x) + O_{L^\infty}(|x|^{-m_j - 1}), 
\eeq
where $m_j > 3$ and 
\beq
V_0^{(j)}(x) = |x|^{-m_j} v_0^{(j)}(\theta) \text{ with } v_0^{(j)} \neq 0 \text{ for } |x|>1. 
\eeq 
We divide the proof into two steps. 

{\bf Step 1:} We claim that $m_1 = m_2$. We argue by contradiction and assume without loss of generality that $m_1 < m_2$. Then 
\beq
V^{(1)}(x) - V^{(2)}(x) = V_0^{(1)}(x) + O_{L^\infty}( |x|^{-m_1-\beta})
\eeq
where $\beta = \min(1, m_2 - m_1)>0$.  
We can apply Lemma \ref{lm-uasymp} and \ref{lm-uasymp1}  to get that for $|x| >0$ sufficiently large, 
\beqq\label{eq-r0-1}
R_0^+(\la)((V^{(1)} - V^{(2)})\Phi_0)(|x|\theta) =    |x|^{-1} e^{i \la |x|} (f_0(\la, \theta, w)  + f_{0, 1}(\la, \theta, w)) + L^2(\mbr^3), 
\eeqq
where 
\begin{enumerate}
\item $f_0(\la, \theta, w) = (-i) (2\pi)^{-2} \la  \hat V_0^{(1)}(\la(\theta -  w))$. For fixed $\theta,$ $f_0$ regarded as a function on $\mbs^2$ does not belong to $H^{m_1-2}(\mbs^2)$ unless $XV_0^{(1)}(x, \theta)= 0$ for $x\in \theta^\perp$ and $|x|>1$. 
\item $f_{0, 1}(\la, \theta, w)$ belongs to $H^{m_1-2}(\mbs^2)$ for fixed $\theta.$ 
\end{enumerate}
Similarly, note that $w^{(1)} \in O_{L^\infty}(|x|^{-1})$ and $V^{(1)} - V^{(2)} \in O_{L^\infty}(|x|^{-{m_1}})$. We can use  Lemma \ref{lm-uasymp1} to get 
\beqq\label{eq-r0-2}
R_0^+(\la)((V^{(1)} - V^{(2)})w^{(1)})(|x|\theta) = |x|^{-1} e^{i \la |x|} f_1(\la, \theta, w)   + L^2(\mbr^3), 
\eeqq
where $f_1$ belongs to $H^{m_1-2}(\mbs^2)$ for fixed $\theta.$ Also, 
\beqq\label{eq-r0-3}
R_0^+(\la)( V^{(2)}(w^{(1)} - w^{(2)})) (|x|\theta) = |x|^{-1} e^{i \la |x|} f_2(\la, \theta, w)   + L^2(\mbr^3), 
\eeqq
where $f_2$ belongs to $H^{m_2 -2}(\mbs^2) \subset H^{m_1 -2}(\mbs^2)$ for fixed $\theta.$ 

Now the order $|x|^{-1}$ terms in \eqref{eq-mainid} must vanish in view of \eqref{eq-r0-1}-\eqref{eq-r0-3} because the other terms in \eqref{eq-mainid} are either $L^2(\mbr^3)$ or $o_{L^\infty}(|x|^{-1})$. So we get $f_0 + f_{0, 1} + f_1 + f_2 = 0$. Thus, $f_0$ must belong to $H^{m-2}(\mbs^2)$ for any $\theta\in \mbs^2$ so 
\beqq\label{eq-ray0}
XV_0^{(1)}(x, \theta) = 0
\eeqq
for any $\theta\in \mbs^2, x\in \theta^\perp$ and $|x|>1$. This gives the X-ray transform of $V_0^{(1)}$ along any lines that do not intersect $|x|<1$. In particular, for any plane $\mcp$ in $\mbr^3$ not intersecting $|x|<1$, we see that the X-ray transform of $V_0^{(1)}$ on $\mcp$ is identically zero. Note that $V_0^{(1)}$ is integrable when restricted to $\mcp$. Thus, we can apply the injectivity for X-ray transform on integrable functions (see \cite{Sol}) to conclude that $V_0^{(1)} = 0$ on $\mcp$. By varying $\mcp,$ we see that $V_0^{(1)}(x) = 0$ for $|x|>1$. This contradicts to the assumption $V_0^{(1)}(x) \neq 0$ for $|x|>1.$ Thus, $m_1 = m_2$. 


{\bf Step 2:} Now we can assume $m_1 = m_2 = m$. Then 
 \beq
 \begin{split}
 V^{(1)}(x) - V^{(2)}(x) &= V_0^{(1)}(x) - V_0^{(2)}(x) + O_{L^\infty}(|x|^{-m-1})\\
  &=  |x|^{-m} (v_0^{(1)}(\theta) - v_0^{(2)}(\theta)) + O_{L^\infty}(|x|^{-m-1}), \quad |x|>1. 
  \end{split}
\eeq 
We consider the identity \eqref{eq-mainid} and repeat the argument in  Step 1. As in \eqref{eq-ray0}, we can conclude that
\beq 
X(V_0^{(1)} - V_0^{(2)})(x, \theta) = 0
\eeq
for any $\theta\in \mbs^2, x\in \theta^\perp$ and $|x|>1$. By the same argument in Step 1, we get $V_0^{(1)}(x) = V_0^{(2)}(x)$ for $|x|>1$ so $V^{(j)} \in S^{m+1}_{\phg}(\mbr^3)$. We continue by induction to conclude that $V^{(1)}- V^{(2)} \in S^{-\infty}_{\phg}(\mbr^3)$. This completes the proof of Theorem \ref{thm-main}.

 \section{The boundary pairing}\label{sec-bdy}
 We prove  Lemma \ref{lm-bp2} in this section.  We first recall the boundary pairing for the Laplacian, see for example Proposition 3.1 of \cite{UhVa}. 
 \begin{lemma}
\label{lm-bp1}
Suppose $V\in L^\infty_{\comp}(\mbr^3)$ is real valued and $\la>0$. Suppose $(-\lap + V - \la^2) u_\pm \in L^{2, s}(\mbr^3), s> 1/2$ and 
\beq
\begin{gathered}
u_+ = e^{-i\la |x|} |x|^{-1}g_{+-}(\theta) + e^{i \la |x|}|x|^{-1}g_{++}(\theta) + L^2(\mbr^3),\\
u_- = e^{-i\la |x|} |x|^{-1}g_{--}(\theta) + e^{i \la |x|}|x|^{-1}g_{-+}(\theta) + L^2(\mbr^3), 
\end{gathered}
\eeq
where $g_{\pm\pm}\in C^\infty(\mbs^2)$. Then
\beqq\label{eq-par} 
\begin{gathered}
\langle u_+, (-\lap + V - \la^2)u_-\rangle -  \langle (-\lap + V - \la^2)u_+, u_-\rangle 
= 2i  \la (\langle g_{++}, g_{-+}\rangle - \langle g_{+-}, g_{--}\rangle).
\end{gathered}
\eeqq  
\end{lemma}
 An essential part of the proof of Lemma \ref{lm-bp1} in \cite{UhVa} (see also \cite{Mel}) is to compute the left hand side of \eqref{eq-par} on $B(R) = \{x\in \mbr^3: |x|< R\}$ for some $R>0$ by using integration by parts. 
 For the fractional Laplacian, integration by parts formula on bounded domains is known in the literature, see for example Theorem 5.1 of \cite{Gru}. However, such formula is valid for functions in the so-called $\mu$-transmission spaces which are usually singular at $\p B(R)$ hence not applicable to our setting. We will prove Lemma \ref{lm-bp2} by reducing it to Lemma \ref{lm-bp1}. 
 
 \bpf[Proof of Lemma \ref{lm-bp2}]
 Assume that  $((-\lap)^\ha - \la) u_\pm = f_\pm \in  L^{2, 1}(\mbr^3)$. Note that $u_\pm\in L^{2, -1}(\mbr^3)$ because of the expression \eqref{eq-upm}. Thus  $\langle ((-\lap)^\ha -\la)  u_+, u_-\rangle$ makes sense as an $L^2$ pairing. For $\eps>0$, we see that  
\beq
\begin{gathered}
\langle ((-\lap)^\ha -\la)  u_+, u_-\rangle =  \lim_{\eps\rightarrow 0+} \langle ((-\lap)^\ha - \la) u_+, e^{-\eps|x|} u_-\rangle  
\end{gathered}
\eeq
by dominated convergence. 
Next, we can verify that $w_\pm = ((-\lap)^\ha + \la)f_\pm$ belongs to $L^{2, 1}(\mbr^3)$. In fact, it suffices to show that $x^j w_\pm \in L^2(\mbr^3)$ for $j = 1, 2, 3.$ But this follows from the fact that $\p_{\xi_j}\msf(f_\pm)(\xi)\in L^2(\mbr^3)$. Thus we have $(-\lap - \la^2) u_\pm = w_\pm \in L^{2, 1}(\mbr^3)$. Note that $((-\lap)^\ha + \la)^{-1}: L^2(\mbr^3)\rightarrow L^2(\mbr^3)$ is bounded for $\la >0$. Thus we can use Plancherel's formula to get 
\beq
\begin{split}
\langle ((-\lap)^\ha -\la)  u_+, u_-\rangle   &= \lim_{\eps\rightarrow 0+}   \langle (-\lap - \la^2) u_+, ((-\lap)^\ha + \la)^{-1} (e^{-\eps|x|} u_-)\rangle\\
 &=  \langle (-\lap - \la^2) u_+, ((-\lap)^\ha + \la)^{-1} u_-\rangle. 
\end{split}
\eeq 
For the last equality, we need to show $((-\lap)^\ha + \la)^{-1}u_-\in L^{2, -1}(\mbr^3)$ and use dominated convergence. In fact, we will show below that $((-\lap)^\ha + \la)^{-1} u_-$ has an asymptotic expansion of the form
\beqq\label{eq-vasymp}
\begin{gathered}
 -(2\la)^{-1}e^{-i\la |x|}  |x|^{-1} g_{--} (\theta)  - (2\la)^{-1} e^{i \la |x|} |x|^{-1} g_{-+}(\theta) + L^2(\mbr^3).
 \end{gathered}
\eeqq
 Then we can apply Lemma \ref{lm-bp1} to get 
 \beq
\begin{split}
&\langle ((-\lap)^\ha -\la)  u_+, u_-\rangle   \\
  =  &\langle  u_+, (-\lap - \la^2)  ((-\lap)^\ha + \la)^{-1} u_-\rangle  
- (2\la)^{-1} 2i  \la (\langle g_{++}, g_{-+}  \rangle - \langle g_{+-}, g_{--} \rangle)\\
  = & \langle u_+,  ((-\lap)^\ha - \la) u_-\rangle - i (\langle g_{++}, g_{-+}  \rangle - \langle g_{+-}, g_{--} \rangle), 
\end{split}
\eeq 
which is the desired formula in Lemma \ref{lm-bp2}.

The rest of the proof is devoted to proving \eqref{eq-vasymp}. Assume that  
\beqq\label{eq-v}
\begin{gathered}
v = v_- + v_+ + L^2(\mbr^3), \text{ where } v_\pm(|x|\theta) = e^{\pm i\la |x|} |x|^{-1}h_\pm(\theta)
\end{gathered}
\eeqq
and $h_\pm \in C^\infty(\mbs^2).$ We show $((-\lap)^\ha + \la)^{-1}v\in L^{2, -1}(\mbr^3)$ and compute its asymptotics at infinity. We know that $((-\lap)^{\ha} + \la)^{-1} w \in L^2(\mbr^3)$ if $w\in L^2(\mbr^3)$. So it suffices to consider $((-\lap)^\ha + \la)^{-1}v_\pm$. Below we analyze the $+$ case. The $-$ case is similar.

Let $\eps> 0$ and $\phi_\eps(x) = e^{-\eps|x|}$. We  compute $((-\lap)^\ha + \la)^{-1}(\phi_\eps v_+)$ for $\eps>0.$ In the end, we consider the limit as $\eps \rightarrow 0+$ in the sense of distributions.  Let $R>0$. We use polar coordinates $\xi = \tau \eta, \tau >0, \eta\in \mbs^2$ and $x = \rho w, \rho>0, w\in \mbs^2$ to obtain that 
\beq
\begin{split}
\widehat{\phi_\eps v_+}(\tau \eta) & = \int_{\mbs^2}\int_0^\infty e^{-i \rho \tau w\cdot \eta} e^{ i\la  \rho}e^{-\eps \rho}h_+(w)  \rho  d\rho dw \\
&= \int_{\mbs^2}\int_0^R e^{-i \rho \tau w\cdot \eta} e^{ i\la  \rho}e^{-\eps \rho}h_+(w)  \rho  d\rho dw + \int_{\mbs^2}\int_R^\infty e^{-i \rho \tau w\cdot \eta} e^{ i\la  \rho}e^{-\eps \rho}h_+(w)  \rho  d\rho dw \\
& \doteq A_+(\tau \eta) + B_+(\tau\eta).  
\end{split}
\eeq
Note that $A_+\in L^2(\mbr^3)$ for all $\eps\geq 0$. We focus on $B_+(\tau\eta)$ below. We consider the behavior of $B_+(\tau\eta)$ for $\tau\in (0, 1)$ and $\tau \in (1, \infty)$. 

If $\tau >1$, we can use the stationary phase lemma for the integration in $w$ to get for $\rho \rightarrow \infty$  that 
\beqq\label{eq-temp0}
\begin{split}
B_+(\tau \eta) 
 &= \int_{\mbs^2}\int_{\tau R}^\infty e^{-i \rho  w\cdot \eta} e^{ i\la  \rho/\tau }e^{-\eps \rho/\tau }h_+(w)  \tau^{-2} \rho  d\rho dw\\
&= (-2\pi i)  \int_{\tau R}^\infty \rho^{-1} e^{i \rho}e^{-\eps \rho/\tau} e^{ i\la  \rho/\tau} \sum_{j = 0}^\infty \rho^{-j} h_{+, j}(-\eta)   \tau^{-2} \rho  d\rho \\
&+ (2\pi i)  \int_{\tau R}^\infty \rho^{-1} e^{-i \rho} e^{ i\la  \rho/\tau}e^{-\eps \rho/\tau}   \sum_{j = 0}^\infty \rho^{-j}  h_{+, j}(\eta)  \tau^{-2} \rho  d\rho + L^2(\mbr^3),
 \end{split}
\eeqq
where $h_{+, j}\in C^\infty(\mbs^2)$ and $h_{+, 0} = h_{+}$. Here, we added an $L^2(\mbr^3)$ term because the asymptotic expansion holds for $\rho> R'$ for some $R'>0.$ But for $\tau < R'/R$ and $\tau>1$, $B_+$ is only $L^2(\mbr^3)$. We observe that for $j\geq 1$, the terms in the two summations of \eqref{eq-temp0} yield $L^2$ functions in $(\tau, \eta)$ (for $\tau >1$), uniformly for $\eps \geq 0$. In fact, for $j\geq 2$, we estimate that 
\beqq\label{eq-temp}
\begin{gathered}
|\int_{\tau R}^\infty \rho^{-1} e^{\pm i \rho}e^{-\eps \rho/\tau} e^{ i\la  \rho/\tau}   \rho^{-j} h_{+, j}(\mp \eta)   \tau^{-2} \rho  d\rho |
\leq C \tau^{-2} |\int_{\tau R}^\infty  \rho^{-j}    d\rho | \leq C \tau^{-2-j+1}
\end{gathered}
\eeqq
and we note that $\tau^{-1-j}$ as a function of $(\tau, \eta)$ is $L^2$ for $j\geq 2.$  For $j = 1,$ the integral in $\rho$ in the middle term of \eqref{eq-temp} is not finite. However, we note that  
\beq
a(\tau) = \int_R^\infty \rho^{-1} e^{i(\tau -\la) \rho} d\rho
\eeq
is an $L^2$ function for $\tau \in \mbr$. Thus $\tau^{-2}a(\tau)$ is $L^2$ in $(\tau, \eta)$ variables. Thus we proved the claim. 
We can write 
\beqq\label{eq-Bp}
\begin{split}
B_+(\tau \eta) &= (-2\pi i)  \int_{\tau R}^\infty  \tau^{-2} e^{i \rho } e^{-\eps \rho/\tau}e^{ i\la  \rho/\tau}  h_{+}(-\eta)   d\rho \\
&+ (2\pi i)  \int_{\tau R}^\infty  \tau^{-2} e^{-i \rho} e^{ i\la\rho/\tau }   e^{-\eps \rho/\tau} h_{+}(\eta)    d\rho + L^2(\mbr^3)
 \end{split}
\eeqq 
for $\tau >1$. We further get that 
\beq
\begin{split}
(-2\pi i)  \int_{\tau R}^\infty  \tau^{-2} e^{i \rho } e^{-\eps \rho/\tau}e^{ i\la  \rho/\tau}  h_{+}(-\eta)   d\rho 
&= (-2\pi i)  \int_{R}^\infty  \tau^{-1} e^{i \rho\tau } e^{-\eps \rho }e^{ i\la  \rho }  h_{+}(-\eta)   d\rho  \\
&= 2\pi \frac{e^{iR(\tau + \la + i\eps)}}{\tau (\tau + \la + i\eps)} h_+(-\eta).
\end{split}
\eeq 
The other integral in \eqref{eq-Bp} is similar. Thus we conclude that for $\tau >1, \eta\in \mbs^2$, 
\beqq\label{eq-Bfin}
\begin{gathered}
B_+(\tau \eta) = I_1 + I_2 + L^2(\mbr^3) \text{ where }\\
I_1 =  2\pi \frac{e^{iR(\tau + \la + i\eps)}}{\tau (\tau + \la + i\eps)} h_+(-\eta), \quad I_2 = -2\pi \frac{e^{iR(-\tau + \la + i\eps)}}{\tau (-\tau + \la + i\eps)} h_+(\eta). 
\end{gathered}
\eeqq 

Now consider $B_+(\tau \eta)$ for $\tau \in (0, 1)$. We split the integral in $B_+$ as 
\beq
\begin{gathered}
B_+(\tau \eta) 
 =  \int_{\mbs^2}\int_R^{R/\tau} e^{-i \rho \tau w\cdot \eta} e^{ i\la  \rho}e^{-\eps \rho}h_+(w)  \rho  d\rho dw
 +  \int_{\mbs^2}\int_{R/\tau}^\infty e^{-i \rho \tau w\cdot \eta} e^{ i\la  \rho}e^{-\eps \rho}h_+(w)  \rho  d\rho dw. 
 \end{gathered}
\eeq
The computation of the second integral is  the same as for the case of $\tau>1$ and we get the same expression as \eqref{eq-Bfin}. For the first integral, we get 
\beq
\begin{gathered} 
 \int_{\mbs^2}\int_R^{R/\tau} e^{-i \rho \tau w\cdot \eta} e^{ i\la  \rho}e^{-\eps \rho}h_+(w)  \rho  d\rho dw 
 = \int_{\mbs^2} [ \frac{1}{i (-\tau w\cdot \eta + \la +i\eps )} \rho e^{i\rho (-\tau w\cdot \eta + \la +i\eps)}|_{R}^{R/\tau}  \\
-  \frac{1}{i (-\tau w\cdot \eta + \la +i\eps )} \int_R^{R/\tau}  e^{i\rho (-\tau w\cdot \eta + \la +i\eps)} d\rho] dw=    O(1/\tau)
 \end{gathered}
\eeq
as $\tau \rightarrow 0$ for all $\eps \geq 0$. Note that the singularity is $L^2(\mbr^3)$ integrable. Thus we conclude that \eqref{eq-Bfin} holds for all $\tau >0.$ 
 
Finally,  we compute in polar coordinates $x = r\theta, \xi = \tau \eta$ where $r, \tau >0, \theta, \eta\in \mbs^2$ that 
\beq
\begin{split}
((-\lap)^{\ha}+ \la)^{-1} (\phi_\eps v)(r\theta) 
&= (2\pi)^{-3}\int_{\mbr^n} e^{ir\theta\cdot \xi} \frac{1}{|\xi|+ \la} \widehat {\phi_\eps v_+}(\xi) d\xi \\
& = (2\pi)^{-3} \int_{\mbs^2}  \int_{0}^\infty e^{i r \tau \theta \cdot \eta}\widehat {\phi_\eps v_+}(\tau \eta)  \frac{ \tau^2}{\tau + \la}  d\tau  d\eta.   
\end{split}
\eeq
We use the expression \eqref{eq-Bfin} and compute the above for $I_1, I_2$. It is easy to see that $I_1$ is in $L^2(\mbr^3)$ for $\eps \geq 0$. For the second term, we have 
\beq
\begin{split}
((-\lap)^{\ha}+ \la)^{-1}  I_2(r\theta) &= -(2\pi)^{-2} \int_{\mbs^2}  \int_{0}^\infty e^{i r \tau \theta \cdot \eta}  \frac{e^{iR(-\tau + \la + i\eps)}}{\tau (-\tau + \la + i\eps)} h_+(\eta) \frac{ \tau^2}{\tau + \la}  d\tau  d\eta  \\
 &= -(2\pi)^{-2} \int_{\mbs^2}  \int_{0}^\infty e^{i r \tau \theta \cdot \eta}  \frac{\tau e^{iR(-\tau + \la + i\eps)}}{(\tau + \la)(-\tau + \la + i\eps)} h_+(\eta)   d\tau  d\eta.
\end{split}
\eeq
Now we use the same argument as in Lemma \ref{lm-uasymp}. Consider a cut-off function $\psi_0$ on $\mbs^2$ such that $\psi_0(w)$ is identically equal to $1$ in a neighborhood of $w\cdot \theta = 0$, and $\psi_0(w)$ vanishing identically in neighborhoods of both $w = \pm \theta$. Let $\psi_0, \psi_\pm$ be a partition of unity with $\psi_\pm$ supported in $\pm \theta \cdot w \geq 0$ such that $\psi_0 + \psi_++\psi_- = 1$ on $\mbs^2$. Then we write 
\beq
\begin{gathered}
((-\lap)^{\ha}+ \la)^{-1}  I_2(r\theta) = \tilde I_{2, 0}(r\theta) + \tilde I_{2, +}(r\theta) + \tilde I_{2, -}(r\theta), \text{ where }\\
\tilde I_{2, t}(r\theta) = -(2\pi)^{-2} \int_{\mbs^2}  \int_{0}^\infty e^{i r \tau \theta \cdot \eta}  \frac{\tau e^{iR(-\tau + \la + i\eps)}}{(\tau + \la)(-\tau + \la + i\eps)} h_+(\eta)   \psi_t(\eta)d\tau  d\eta. 
\end{gathered}
\eeq
We observe that $\tilde I_{2, 0}, \tilde I_{2, -}$ are $L^2$ uniformly for $\eps \geq 0.$ For $\tilde I_{2, +}$, we can use again the contour deformation in Lemma \ref{lm-uasymp} to compute that 
\beq
\begin{gathered}
\tilde I_{2, +}(r\theta)  
 =  -(2\pi)^{-2} (2\pi i) \int_{\mbs^2}   e^{i r  (\la + i\eps) \theta \cdot \eta}  \frac{\la + i\eps}{2\la +i\eps} h_+(\eta)   \psi_+(\eta)   d\eta + L^2(\mbr^3). 
 \end{gathered}
\eeq
Note that our computation are uniform for $\eps >0$ so far. We can take the limit as $\eps \rightarrow 0+$ in the sense of distributions. We finally use the stationary phase lemma to get 
\beq
\begin{gathered}
\tilde I_{2, +}(r\theta) 
  = -(2\pi)^{-2} (2\pi i) (-2\pi i)  \frac{1}{2r\la} e^{i r  \la}  h_+(\theta)    + L^2(\mbr^3)
   = -\frac{1}{2\la r} e^{i r  \la }   h_+(\theta)    + L^2(\mbr^3). 
 \end{gathered}
\eeq  
To summarize,  we proved that 
\beq
\begin{gathered}
((-\lap)^{\ha}+ \la)^{-1} v_+(r\theta)  =   -e^{i\la r} r^{-1}  h_{+}(\theta) (2\la)^{-1}   + L^2(\mbr^3). 
\end{gathered}
\eeq
This completes the calculation for $((-\lap)^{\ha}+ \la)^{-1} v_+$. Repeating the same argument to $v_-$, we get that 
\beq
\begin{gathered}
((-\lap)^{\ha}+ \la)^{-1} v_-(r\theta)   =    -e^{-i\la r} r^{-1}  h_{-}(\theta) (2\la)^{-1}   + L^2(\mbr^3). 
\end{gathered}
\eeq
This finishes the proof of \eqref{eq-vasymp}. 
\epf

\section*{Acknowledgement}
The authors wish to thank Katya Krupchyk for helpful conversations and providing reference \cite{Gru}. GU and YW are both supported by NSF. 


\end{document}